\documentclass{article}
\usepackage{amssymb}
\usepackage{amsbsy}
\begin{document}
\centerline{\LARGE \bf } \vskip 6pt

\begin{center}{\Large \bf The mixed affine quermassintegrals}\footnote{Research is supported by
National Natural Sciences Foundation of China (10971205, 11371334).}\end{center}
\vskip 10pt

\centerline{Chang-Jian Zhao}

\begin{center}{\it Department of Mathematics,
China Jiliang University, Hangzhou 310018, P. R.
China}\end{center} \centerline{\it Email: chjzhao@163.com}

\begin{center}
\begin{minipage}{12cm}
{\bf Abstract} In this paper, we introduce first the mixed affine quermassintegrals of $j$ convex bodies. The Aleksandrov-Fenchel inequality for the mixed affine quermassintegrals of $j$ convex bodies is established. As a application, the Minkowski's, Brunn's Minkowski's inequalities for the mixed affine quermassintegrals are also derived.


\end{minipage}
\end{center}
\vskip 20pt

\noindent{\large \bf 1.~ Introduction}\vskip 10pt

Lutwak [1] proposed to define the affine quermassintegrals for a
convex body $K$, $\Phi_{0}(K)$, $\Phi_{1}(K)$, $\ldots,$
$\Phi_{n}(K)$, by taking $\Phi_{0}(K):=V(K),
\Phi_{n}(K):=\omega_{n}$ and for $0<j<n$,
$$\Phi_{n-j}(K):=\omega_{n}\left[\int_{G_{n,j}}\left(\frac{{\rm
vol}_{j}(K|\xi)}{\omega_{j}}\right)^{-n}
d\mu_{j}(\xi)\right]^{-1/n},\eqno(1.1)$$ where $G_{n,j}$ denotes
the Grassman manifold of $j$-dimensional subspaces in ${\Bbb
R}^{n}$, and $\mu_{j}$ denotes the gauge Haar measure on
$G_{n,j}$, and ${\rm vol}_{j}(K|\xi)$ denotes the $j$-dimensional
volume of the positive projection of $K$ on $j$-dimensional
subspace $\xi\subset {\Bbb R}^{n}$ and $\omega_{j}$ denotes the
volume of $j$-dimensional unit ball (see [2]). Lutwak showed the Brunn-Minkowski inequality for the affine
quermassintegrals. If $K$ and $L$ are convex bodies and $0<j<n$,
then
$$\Phi_{n-j}(K+L)^{1/j}\geq\Phi_{n-j}(K)^{1/j}+\Phi_{n-j}(L)^{1/j}.\eqno(1.2)$$

In this paper, we introduce the mixed affine quermassintegrals of $j$ convex bodies. The Aleksandrov-Fenchel inequality for the mixed affine quermassintegrals of $j$ convex bodies is established. As a application, and the Minkowski inequality is also derived.

\vskip 10pt \noindent{\large \bf 2.~ The mixed affine quermassintegrals}\vskip 10pt

In the section, we introduce first the following concept and  List its properties and related inequalities.

{\bf Definition 2.1} (The mixed affine quermassintegrals of $j$ convex bodies) The mixed
affine quermassintegral of $j$ convex bodies $K_{1},\ldots,K_{j}$,
denoted by $\Phi_{n-j}(K_{1},\ldots,K_{j})$, defined by
$$\Phi_{n-j}(K_{1},\ldots,K_{j}):=\omega_{n}\left[\int_{G_{n,j}}\left(\frac{{\rm
vol}_{j}((K_{1},\ldots,K_{j})|\xi)}{\omega_{j}}\right)^{-n}
d\mu_{j}(\xi)\right]^{-1/n},\eqno(2.1)$$ where $0\leq j\leq n$.

When $K_{1}=\cdots=K_{j}=K$, $\Phi_{n-j}(K_{1},\ldots,K_{j})$ becomes Lutwak's affine quermassintegral $\Phi_{n-j}(K)$. When $K_{1}=\cdots=K_{j-1}=K$ and $K_{j}=L$, $\Phi_{n-j}(K_{1},\ldots,K_{j})$ becomes a new affine geometric quantity, denoted by $\Phi_{n-j}(K,L)$ and call it mixed affine quermassintegral of $K$ and $L$. When $K_{1}=\cdots=K_{j-i-1}=K$, $K_{j-i+1}=\cdots=K_{j-1}=B$ and $K_{j-i}=L$, $\Phi_{n-j}(K_{1},\ldots,K_{j})$ becomes another new affine geometric quantity, denoted by $\Phi_{n-j,i}(K,L)$ and call it $i$-th mixed affine quermassintegral of $K$ and $L$, where $0\leq i<j\leq n$. When $K_{1}=\cdots=K_{j-i}=K$ and $K_{j-i+1}=\cdots=K_{j}=B$, $\Phi_{n-j}(K_{1},\ldots,K_{j})$ becomes a new affine geometric quantity, denoted by $\Phi_{n-j,i}(K)$ and call it $i$-th mixed affine quermassintegral of $K$, where $0\leq i<j\leq n$.

Obviously, the mixed affine quermassintegrals of $j$ convex bodies is invariant under simultaneous unimodular centro-affine transformation.

{\bf Lemma 2.1} {\it If
$K_{1},\ldots,K_{j}\in {\cal K}_{o}^{n}$ and $0\leq j\leq n$, then for any $g\in{\rm SL(n)}$}
$$\Phi_{n-j}(gK_{1},\ldots,gK_{j})=\Phi_{n-j}(K_{1},\ldots,K_{j}).$$

As we all know, according to the Brunn-Minkowski theory, a very natural question is raised: are there some isoperimetric inequalities about the mixed affine quermassintegrals of $j$ convex bodies? The following perfectly answers the question and establish Minkowski's, and leksandrov-Fenchel's and Brunn-Minkowski's inequalities for the mixed affine quermassintegrals.

{\bf Theorem 2.1} {\rm (The Minkowski inequality for mixed affine quermassintegrals)} {\it If $K,L\in {\cal K}^{n}_{o}$ and $0\leq j\leq n$, then
$$\Phi_{n-j}(K,L)^{j}\geq\Phi_{n-j}(K)^{j-1}\Phi_{n-j}(L),\eqno(2.2)$$
with equality if and only if $K$ and $L$ are homothetic.}

{\it Proof}~ This follows immediately from the Minkoweski's, and H\"{o}lder's inequalities.\hfill$\Box$

Next, we establish an Aleksandrov-Fenchel inequality for the mixed affine quermassintegral of $j$ convex bodies $K_{1},\cdots,K_{j}$.

{\bf Theorem 2.3} {\rm (The Aleksandrov-Fenchel inequality for mixed affine quermassintegrals of $j$ convex bodies)} {\it If $K_{1},\cdots,K_{j}\in {\cal K}^{n}_{o}$, $0\leq j\leq n$ and $0<r\leq j$, then}
$$\Phi_{n-j}(K_{1},\cdots,K_{j})\geq
\prod_{i=1}^{r}\Phi_{n-j}(K_{i},\cdots,K_{i},K_{r+1},\cdots,K_{j})^{1/r}.\eqno(2.3)$$

{\it Proof}~ This follows immediately from the Aleksandrov-Fenchel inequality and H\"{o}lder's inequality.\hfill$\Box$

Unfortunately, the equality conditions of the Aleksandrov-Fenchel inequality are, in general, unknown.

{\bf Corollary 2.1} {\it If $K_{1},\cdots,K_{j}\in {\cal K}^{n}_{o}$ and $0\leq j\leq n$, then
$$\Phi_{n-j}(K_{1},\cdots,K_{j})^{j}\geq\Phi_{n-j}(K_{1})\cdots\Phi_{n-j}(K_{j}),\eqno(2.4)$$
with equality if and only if $K_{1},\cdots,K_{j}$ are homothetic.}

{\bf Proof}
The special case $r=j-1$, of inequality (4.3), is
$$\Phi_{n-j}(K_{1},\cdots,K_{j})^{j-1}\geq\Phi_{n-j}(K_{1},K_{j})\cdots\Phi_{n-j}(K_{j-1},K_{j}).$$
When above inequality is combined with the Minkowski inequality (4.2), the
result is
$$\Phi_{n-j}(K_{1},\cdots,K_{j})^{j}\geq\Phi_{n-j}(K_{1})\cdots\Phi_{n-j}(K_{j}),$$
with equality if and only if $K_{1},\cdots,K_{j}$ are homothetic.
\hfill$\Box$

Finally, we simply prove the Brunn-Minkowski inequality for the affine quermassintegrals by using the mixed affine quermassintegrals theory introduced in this section.

{\bf Theorem 2.3} {\rm (The Brunn-Minkowski inequality for affine quermassintegrals)} {\it If $K,L\in {\cal K}^{n}_{o}$ and $0\leq j\leq n$, then for} $\varepsilon>0$
$$\Phi_{n-j}(K+\varepsilon\cdot L)^{1/j}\geq\Phi_{n-j}(K)^{1/j}+\varepsilon\Phi_{n-j}(L)^{1/j}.\eqno(2.5)$$

{\it Proof}~ This follows immediately from (2.1) and (2.2).\hfill$\Box$

\end{document}